\input amstex
\documentstyle{amsppt}
\magnification=\magstep1 \pagewidth{6.2 in} \pageheight{7.7in}
\hcorrection{-0.4in} \vcorrection{-0.4in} \abovedisplayskip=10pt
\belowdisplayskip=10pt
\parskip=4pt
\parindent=5mm
\baselineskip=2pt
\NoBlackBoxes

\topmatter

\title A note on the Fourier transform of $p$-adic $q$-integrals
\endtitle
\author   Taekyun Kim\endauthor
\affil{{\it Department of Mathematics Education}\\
{\it Kongju University, Kongju 314-701, Korea }\\
{\it e-mail:tkim$\@$kongju.ac.kr}}\\
\endaffil

\define\Q{\Bbb Q_p}
\define\C{\Bbb C_p}
\define\BZ{\Bbb Z}
\define\BQ{\Bbb Q}
\define\BC{\Bbb C}
\define\BN{\Bbb N}

\define\Z{\Bbb Z_p}

 \subjclass 11B68, 11S80
\endsubjclass \abstract{The $p$-adic $q$-integral (=
$I_q$-integral) was defined by author in the previous paper [1,
3].
 In this paper, we consider $I_q$-Fourier transform and investigate some properties
 which are related to this transform.}
\endabstract
\rightheadtext{ T. Kim} \leftheadtext{A note on $I_q$-Fourier
transform } \TagsOnRight
\endtopmatter
\document
\head \S 1. Introduction \endhead
 Let us denote $\BN, \BZ,\BQ, \BC$ sets of positive integer, integer,
rational and complex numbers respectively. Let $p$ be  prime and
$x \in \BQ$. Then $x= p^{v(x)} {m \over n}$, where $m,n, v=v(x)
\in \BZ$, $m$ and $n$ are not divisible by $p$. Let $|x|_p =
p^{-v(x)}$ and $ |0|_p =0$. Then $|x|_p$ is valuation on $\BQ$
satisfying
$$ |x+y|_p \leq \max \{ |x|_p , |y|_p \} .$$
Completion of $\BQ$ with respect to $| \cdot |_p$ is denoted by
$\Q$ and called the field of $p$-adic rational numbers. $\C$ is
the completion of algebraic closure of $\Q$ and $\Z = \{ x \in \Q
~|~ |x|_p  \leq 1 \} $ is called the ring of $p$-adic rational
integers(see 1, 2, 3, 4]). Let $l$ be a fixed integer and let $p$
be a fixed prime number. We set
$$\aligned
&X_l=\varprojlim_N (\BZ/lp^N \BZ), \text{ and $X_1=\Bbb Z_p , $}\\
&X^*=\bigcup\Sb 0<a<lp\\ (a,p)=1\endSb (a+lp\Z ),\\
&a+lp^N \Z =\{x\in X\mid x\equiv a\pmod{lp^N}\},
\endaligned $$
where $N \in \BN$ and  $a\in \BZ$ lies in $0\leq a<lp^N$, cf. [3].

When one talks of $q$-extension, $q$ is considered in many ways
such as an indeterminate, a complex number $q \in \BC$, or a
$p$-adic number $q \in \C$.  In this paper, we assume that $q \in
\C$ with $|q-1|_p < p^{- {1 \over {p-1}}}$, so that $q^x = \exp (x
\log q)$ for each $x \in X$. We use the notation as $[x] = [x; q]=
{{1-q^x}\over {1-q}}$ for each $x \in X$. Hence $ \lim_{q \to 1}
[x] =x$, cf.[3]. For any positive integer $N,$ we set
$$\mu_q (a+lp^N\Z)=\frac{q^a }{[lp^N ]}, \text{ cf. [3]},$$
and this can be extended to a distribution on $X$. This
distribution yields an integral as follows (see [3]):
$$ I_q (f)=\int_{\Z}  f(x) \,d\mu_q(x)=\int_X  f(x) \,d \mu_q (x) ,$$
where $f \in UD(\Bbb Z_p )=\text{ the space of uniformly
differentiable function  on } \Bbb Z_p $ with values in $\Bbb
C_p$, cf. [3]. Let $ C_{p^n}$ be the cyclic group consisting of
all $p^n$-th roots of unity in $\Bbb C_p $ for any $n\geq 0$ and $
T_p$ be the direct limit of $ C_{p^n}$ with respect to the natural
morphism, hence $ T_p$ ia the union of all $C_{p^n}$ with discrete
topology. $U_p$ denotes the group of all principal units in $\Bbb
C_p .$ For any $ f \in UD( \Bbb Z_p, \Bbb C_p )$, we have an
integral $I_0 (f)$ with respect to the so called invariant measure
$\mu_0 $:
$$I_0(f)=\int_{\Bbb Z_p} f(x)d\mu_0(x)=\lim_{n\rightarrow \infty}
\frac{1}{p^n}\sum_{x=0}^{p^n-1}f(x), \text{ cf. [3, 4] }, $$ and
the Fourier transform $\hat{f_w}=I_0(f\phi_w ),$ where $\phi_w$
denotes a uniformly differentiable function on $\Bbb Z_p$
belonging to $w\in T_p$ defined by $\phi_w(x)=w^x, $ cf. [4]. Now
we introduce the convolution for any $f, g \in UD(\Bbb Z_p , \Bbb
C_p )$ due to Woodcock as follows:
$$ f*g(x)=\sum_{w} \hat{f_w}\hat{g_w}\phi_{w^{-1}}(x), \text{ see
[3, 4]. } $$ As known results, $f*g \in UD(\Bbb Z_p , \Bbb C_p )
,$ and $\widehat{(f*g)_w}= \hat{f_w}\hat{g_w} ,$ ( see [4]). In
this paper, we consider $I_q$-Fourier transform and investigate
some properties which are related to this transform. It is easy to
see that $I_q$-Fourier transform is exactly same $I_0$-Fourier
transform when $q=1$.

\head \S 2. $I_q$-Fourier transform
\endhead

 For any $f \in UD(X)$, the $p$-adic $q$-integral
was defined by
$$ I_q (f)  = \int_X f(x) ~d \mu_q (x)
  = \lim_{N \to \infty} {1 \over {[lp^N ]}} \sum_{0\leq x <
 lp^N } f(x) q^x ,
 \text{ cf. [3] .}$$
  Note that
$$ I_0 (f)  = \lim_{q \to 1} I_q (f)= \int_{\Bbb Z_p } f(x) ~d \mu_1
(x)  = \lim_{N \to \infty} {1 \over {lp^N }} \sum_{0\leq x <
 lp^N } f(x) ,$$ and that
$$ I_0 (f_1 ) = I_0 (f) + f' (0), \text{ where $f^{\prime }(0)=\frac{ d}{dx} f(x)|_{x=0} $ and $f_1 (x) =f(x+1)$}.$$
 Let $T_p = \cup_{n \geq 1} C_{p^n} = \lim_{n \to \infty} \BZ / p^n \BZ$,
where $C_{p^n } =\{ \xi \in X | ~ \xi^{p^n} =1 \}$ is the cyclic
group of order $p^n$, see [1]. For $ \xi \in T_p $, we denote by $
\phi_\xi : \Z \to \C$ the locally constant function $ x \mapsto
\xi^x $. If we take $f(x) = \phi_\xi (x) e^{tx}$, then we have
that $ \int_X e^{tx} \phi_\xi (x) d\mu_q (x) = {t + \log q \over
{q\xi e^t -1}} , \text{ cf. [2]}. $ We now consider $I_q$-Fourier
transform as follows:
$$(\hat{f_w})_q =I_q(\phi_{w} f)=\int_{\Bbb Z_p } \phi_{w}(x) f(x)
d\mu_q(x), \text{ where  $f\in UD(\Bbb Z_p)$, $w \in T_p $ },$$
and its inverse transform is derived by
$$\frac{\log q}{q-1}\lim_{n \rightarrow \infty}\sum_{x=0}^{p^n-1} w^{-x}I_q
(\phi_{w} f)=\frac{\log q}{q-1}\lim_{n\rightarrow
\infty}\sum_{x=0}^{p^n-1}w^{-x}\frac{1}{[p^n]}\sum_{z=0}^{p^n-1}w^{z}f(z)q^z
= f(x)\phi_{q} (x) .$$
 Thus, we obtain  the below proposition.
\proclaim{Proposition  1 } Let $f\in UD(\Bbb Z_p, \Bbb C_p)$. Then
we have the inverse formula of $I_q$-Fourier transform as follows:
$$ f(x)\phi_q(x)=\frac{\log q}{q-1}\sum_{w}\phi_{w^{-1}}(\hat{f_w})_{q}, \text{ where $\sum_{w}
=\lim_{n\rightarrow \infty}\sum_{w\in C_p{^n}} $}. $$
\endproclaim
Remark.  In [4], we note that if $\alpha \in U_p $, then
$\phi_{\alpha }$ is called locally analytic character  and if
$\alpha \in T_p ,$ then $\phi_{\alpha}$ is called locally constant
function. For $f, g \in UD ( \Bbb Z_p ) $, we consider  the
convolution of $f, g$ by
$$f *_q g =\sum_{w} (\widehat{f_{wq^{-1}}})_q
(\widehat{g_{wq^{-1}}})_q \phi_{w^{-1}}.$$ Thus, we note that
$$ f *_q g \in UD(\Bbb Z_p ), \text{ and } (\widehat{(f*_q
g)_{wq^{-1}} } )_q =(\widehat { f_{wq^{-1}}})_q (\widehat {
g_{wq^{-1}}})_q .$$ And we also see that
$$ UD(\Bbb Z_p, \Bbb C_p )/\{ f\in UD( \Bbb Z_p , \Bbb C_p )|
f^{\prime } =0 \} \cong C ( \Bbb Z_p , \Bbb C_p ), $$ where
 $C(\Bbb Z_p , \Bbb C_p ) $ is the space of the continuous function
from $\Bbb Z_p $ to $ \Bbb C_p $. Another convolution $\otimes_q$
is induced from the above convolution $ *_q $ by $ (f *_q
g)^{\prime}=-f^{\prime} \otimes_q g^{\prime}. $ Then, we also see
that $f\otimes_q g \in UD(\Bbb Z_p, \Bbb C_p ). $ From these
definitions, we can derive the below theorem.

\proclaim{Theorem 2} For $f, g \in UD(\Bbb Z_p, \Bbb C_p ), $ we
have

$$f *_q g (z) =\frac{q-1}{\log q} I_q^{(x)}(
f(x)g(z-x)\phi_{q^{-1}}(x)) -( f \otimes_q g^{\prime}(z) ),
$$ where $I_q^{(x)}$ means the integration with respect to the
variable $x$.
\endproclaim

Since $(\widehat{ (f*_qg)_{wq^{-1}}})_q  =
(\widehat{f_{wq^{-1}}})_q (\widehat{g_{wq^{-1}}})_q$, for $ w \in
T_p ,$ we have
$$ \int_{\Bbb Z_p } ( f *_q g(x)) q^{-x} d\mu_q(x)=\int_{\Bbb Z_p
}f(x)q^{-x}d\mu_q(x)\int_{\Bbb Z_p}g(x)q^{-x} d\mu_q(x). $$ From
this, we can derive the below worthwhile and interesting formula:

\proclaim{ Theorem 3} Let $ f, g \in UD(\Bbb Z_p , \Bbb C_p ).$
Then we have
$$ I_q^{(z)}(f\otimes g^{\prime}(z)q^{-z} )=\frac{q-1}{\log q}
I_q^{(z)}
(I_q^{(x)}(f(x)g(z-x)q^{-x})q^{-z})-I_q^{(z)}(f(z)q^{-z})I_q^{(z)}(g(z)q^{-z}).$$
\endproclaim

\enddocument